\documentclass[reqno,centertags,12pt]{amsart}
\usepackage{amsmath,amsthm,amscd,amssymb,amsfonts,amscd,latexsym,verbatim,bm}
\usepackage{setspace}
\usepackage{enumitem}
\usepackage{Tabbing}
\usepackage{fancyhdr}
\usepackage{lastpage}
\usepackage{multicol}
\usepackage{wrapfig}
\usepackage{graphicx}
\usepackage{units}
\usepackage{cancel}
\usepackage{mathtools}
\usepackage{accents}
\usepackage{pgfplots}
\usepackage{wrapfig}
\usepackage{hyperref}
\usepackage{graphicx,epsf,cite}

-\marginparwidth \oddsidemargin -4mm \evensidemargin \oddsidemargin

\newtheorem{theorem}{Theorem}[section]

\newtheorem{lemma}[theorem]{Lemma}

\theoremstyle{definition}

\theoremstyle{remark}
\newtheorem{remark}[theorem]{Remark}

\allowdisplaybreaks
\numberwithin{equation}{section}

\newcommand{\lb}{\label}
\newcommand{\f}{\frac}

\newcommand{\supp}{\text{\rm{supp}}}

\newcommand{\abs}[1]{\left \lvert#1\right\rvert}
\newcommand{\tao}{\theta_1}
\newcommand{\taz}{\theta_0}

\newcommand{\beq}{\begin{equation}}
\newcommand{\eeq}{\end{equation}}

\newcommand{\bal}{\begin{align}}
\newcommand{\eal}{\end{align}}
\newcommand{\bals}{\begin{align*}}
\newcommand{\eals}{\end{align*}}

\newcommand{\eps}{\varepsilon}
\newcommand{\del}{\delta}

\newcommand{\al}{\alpha}
\newcommand{\be}{\beta}
\newcommand{\ga}{\gamma}
\newcommand{\laa}{\lambda}

\newcommand{\til}{\tilde}

\newcommand{\rb}[1]{\left(#1\right)}
\newcommand{\ta}{\theta}
\newcommand{\tta}{\til \theta}
\newcommand{\tB}{B_2}
\newcommand{\tphi}{\phi}
\newcommand{\tbe}{\be}

\newcommand{\epsi}{\psi^\eps}
\newcommand{\ela}{\lambda^\eps }
\newcommand{\tl}{{\til \ell}}
\newcommand{\cb}[1]{\left\{#1\right\}}

\newcommand{\bbN}{{\mathbb{N}}}
\newcommand{\bbR}{{\mathbb{R}}}

\newcommand{\mR}{\mathbb{R}}

\title[Transition fronts for ignition-monostable reactions]{Existence and non-existence of transition fronts in mixed ignition-monostable media}
\author{Cole Graham, Tau Shean Lim, Andrew Ma, David Weber}

\address{
Cole Graham: \href{mailto:grahamca@stanford.edu}{grahamca@stanford.edu}\\
Tau Shean Lim: \href{mailto:taushean@math.wisc.edu}{taushean@math.wisc.edu}\\
Andrew Ma: \href{mailto:andyma@math.utexas.edu}{andyma@math.utexas.edu}\\
David Weber: \href{mailto:dsweber@math.ucdavis.edu}{dsweber@math.ucdavis.edu}}

\begin{document}

\begin{abstract}
We study transition fronts for one-dimensional reaction-diffusion equations with compactly perturbed ignition-monostable reactions.
We establish an almost sharp condition on reactions which characterizes the existence and non-existence of fronts.
In particular, we prove that a strong inhomogeneity in the reaction prevents formation of transition fronts, while a weak inhomogeneity gives rise to a front.
Our work extends the results and methods introduced in \cite{Zlatos}, which studies the same question in inhomogeneous KPP media.
\end{abstract}

\maketitle

\section{Introduction}
\lb{intro}

We study transition fronts for one-dimensional reaction-diffusion equations
with \emph{compactly-perturbed ignition-monostable reactions}.
Consider the evolution PDE
\begin{equation}  \label{eq:main}
u_t = u_{xx} + f(x,u),\quad (t,x)\in \mR\times \mR,
\end{equation}
where the nonlinearity $f$ satisfies the following on $\mR\times [0,1]$:
\begin{enumerate}[label=(F\arabic*)]
  \setlength\itemindent{10pt}
\item $f\ge 0$ is Lipschitz continuous with $\ga:=\mbox{Lip}(f)$, and  $f(x,0)=f(x,1) = 0$ for all $x\in \mR$;
  \label{item:general}
  
\item there exists $L>0$ such that $f(x,u)\equiv f_0(u)$ for all $|x|\ge L$, where $f_0$ is an ignition reaction with $f_0\equiv 0$ on $[0,\taz]\cup \{1\}$, $f_0>0$ on $(\taz ,1)$, and $f_0$ is non-increasing on $[1-\tao,1]$ for some $\taz,\tao\in (0,1)$;
  \label{item:ignition}
  
\item \label{item:regularity} the (right hand) derivative $a(x):=f_u(x,0)\ge 0$ exists, and for all $\eps>0$, there exists $\zeta=\zeta(\eps)\in(0,\taz)$ such that
  \begin{equation*}
    (1-\eps)a(x)u\leq f(x,u)\leq (a(x)+\eps)u\quad \mbox{for }(x,u)\in \mR\times [0,\zeta].
  \end{equation*}
\end{enumerate}
As described above, $f$ is obtained by perturbing a homogeneous ignition reaction $f_0$ locally on the interval $[-L,L]$ with an inhomogeneous monostable reaction.
In the present work, we are interested in how such perturbation affects the existence of transition fronts.

The PDE \eqref{eq:main} and its variations are widely used to model a host of natural processes, including thermal, chemical, and ecological dynamics.
By \ref{item:general}, $u\equiv 0, 1$ are two equilibrium solutions of \eqref{eq:main}.
Therefore, one is usually interested in the transition from the (unstable) state $u\equiv 0$ to the (stable) state $u\equiv 1$.
\emph{Transition fronts} are a class of solutions that model this phenomenon.
They are global-in-time solutions $u:\mR^2 \to (0,1)$ of \eqref{eq:main} satisfying
\begin{equation}  \label{item:trans_lims}
\lim_{x\to-\infty} u(t,x)=1, \quad\lim_{x\to+\infty} u(t,x)=0\quad \mbox{for all }t\in \mR
\end{equation}
and the bounded front width condition, that is, for all $\mu \in (0,\f 12 )$,
\begin{equation}  \label{item:trans_width}
\sup _{t\in \mR} L_\mu(t):=\sup_{t\in\mR}\, \mbox{diam} \{x\in\mR|\; \mu\le u(t,x)\le 1-\mu\}<\infty.
\end{equation}
This definition was introduced in \cite{BH12, Matano, Sen}.

The study of transition fronts has seen much activity since the seminal works by Fisher \cite{F37} and Kolmogorov, Petrovskii, and Piskunov \cite{KPP37}, who first studied \emph{traveling fronts} for \eqref{eq:main} with \emph{homogeneous Fisher-KPP reactions}.
Here, traveling fronts are transition fronts of the form $u(t,x)=U(x-ct)$ for some speed $c\in \mR$ and profile $U$ with $\lim_{y\to-\infty}U(y)=1$, $\lim_{y\to\infty}U(y)=0$,
and Fisher-KPP reactions are those $f$ satisfying \ref{item:general}, $f'(0)>0$, and $0<f(u)\le f'(0)u$ for $u\in(0,1)$.
In their work, they found that for each $c\ge c_*:=2\sqrt{f'(0)}$, \eqref{eq:main} admits a unique (modulo translation) traveling front, while no front exists for $c<c_*$.
By means of phase plane analysis, it can be shown that the same existence result holds for general homogeneous monostable reactions ($f(u)>0$ on $(0,1)$), although $c_*\ge 2\sqrt{f'(0)}$ in general (e.g., see \cite{AW78}).
In contrast, for \emph{homogeneous ignition} (defined as in \ref{item:ignition}) and \emph{bistable reactions} (the same as ignition except $f<0$ on $(0,\taz)$ and $\int _0^1 f(u)du>0$), there is only one speed $c_*>0$ which gives rise to a unique (up to translation) traveling front. The unique speed $c_*$ will be called \emph{the spreading speed of }$f$.

Over decades, the study of transition fronts extended to spatially periodic reactions (in which case fronts have time-periodic profiles, and are known as \emph{pulsating fronts}).
Instead of surveying the vast literature, let us refer to the review articles by Berestycki \cite{B03} and Xin \cite{X00}, and the references therein.
The development in general inhomogeneous media is considerably more recent.
The first existence result was obtained by Vakulenko and Volpert \cite{VV11} for small perturbations of homogeneous bistable reactions.
Later, Mellet, Roquejoffre, and Sire \cite{MARS10} proved the existence of fronts for ignition reactions of the form $f(x,u)=a(x)f_0(u)$, where $f_0$ is ignition, and $a(x)$ is bounded with $\inf _{\mR}a(x)>0$, which need not be close to being constant (see also \cite{NR09} for the case of random media, relying on the notion of generalized random traveling waves developed in \cite{Sen}). Zlato\v{s} then extended these results (along with uniqueness and stability) to general inhomogeneous ignition and mixed ignition-bistable media \cite{Z13, ZPreprint}.

Transition fronts has also been investigated in inhomogeneous Fisher-KPP media by several authors. As far as Fisher-KPP reactions are concerned, a strong inhomogeneity in the reaction may prevent existence of transition fronts, while a weak inhomogeneity gives rise to them.
This is translated into the following result proved by Nolen, Roquejoffre, Ryzhik and Zlato\v{s} \cite{Zlatos} for reactions satisfying $0< f(x,u)\le a(x)u$ for all $(x,u)\in \mR\times (0,1)$, with $a(x):=f_u(x,0)$, $a_-:=\inf_{x\in \mR}a(x)>0$, and $a(x)-a_-\in C_c(\mR)$.
They found that when the inhomogeneity of $f$ is strong, in the sense that the principal eigenvalue $\lambda$ of the operator $\partial_{xx}+a(x)$ satisfies $\lambda>2a_-$, any non-constant global-in-time solution $u$ of \eqref{eq:main} is \emph{bump-like} (i.e. $u(t,x)\le C_t e^{-c|x|}$), preventing the existence of transition fronts.
This in fact is the first known example of a reaction function $f$ such that \eqref{eq:main} does not admit any transition front.

Moreover, in the same work, they also show that the existence criterion is (almost) sharp.
In the case of a weak localized inhomogeneity $\lambda<2a_-$, for each $c\in (2\sqrt{a_-},\laa/\sqrt{\laa-a_-})$ the PDE \eqref{eq:main} admits a transition front with \emph{global mean speed} $c$, in the sense that if $X(t):=\sup \{x\in\mR:u(t,x)=\f 12\}$, then
\begin{equation}\lb{gms}
  \lim_{t-s\to\infty}\frac{X(t)-X(s)}{t-s}=c.
\end{equation}
To construct a front, they find an appropriate pair of ordered global-in-time super- and sub-solutions $w\ge v$ that propagate with speed $c$, and recover a front $u$ between them as a locally uniform limit along a subsequence of solutions $(u_n)_{n\in \bbN}$ of the Cauchy problem \eqref{eq:main} with initial data $u_n(-n,\cdot)=w(-n,\cdot)$.
The same method was deployed and extended by Zlato\v{s} \cite{Z12} and by Tao, Zhu, and Zlato\v{s} \cite{TZZ13} to prove the existence of fronts for general inhomogeneous KPP and monostable reactions when $a(x)-a_-$ is not compactly supported.

In the present paper, we modify the approach from \cite{Zlatos} to establish a similar sharp existence criterion for reactions satisfying Hypothesis (F).
As mentioned, such $f$ is obtained by locally perturbing the ignition reaction $f_0$ with a monostable reaction.
We therefore show that a strong perturbation in the reaction prevents the existence of fronts, while a weak perturbation admits them.
The existence criterion in our case is determined by the spreading speed of the reaction $f_0$ and the supremum of the spectrum of the operator $\partial_{xx}+a(x)$.
The spreading speed of $f_0$ is the unique number $c_0>0$ such that the following ODE admits a unique (up to translation) solution:
\begin{equation}
  \lb{tfeq}
  U''+c_0U' +f_0(U)=0,\quad \lim_{x\to -\infty} U(x)=1,\quad \lim_{x\to\infty} U(x)=0.
\end{equation}
On the other hand, the supremum of the spectrum of $\partial_{xx}+a(x)$ is given by
\begin{equation*}
  \lambda := \sup \sigma (\partial_{xx}+a(x))= \sup_{{\psi \in H^1(\mR)}:\,||\psi||_{L^2}=1} {\int_\mR (-[\psi'(x)]^2 + a(x) [\psi (x)]^2)dx}.
\end{equation*}
Since $a(x)\ge 0$ is compactly supported by \ref{item:ignition}, the essential spectrum of $\partial_{xx}+a(x)$ is $(-\infty,0]$, which implies $\lambda\ge 0$.
If $\lambda>0$ (i.e. $a\not\equiv 0$), it is in fact the principal eigenvalue.
Then a corresponding $L^\infty$-normalized principal eigenfunction $\psi$ exists, is unique, and satisfies
\begin{equation}
  \lb{eq:eigen}
  \psi '' + a(x)\psi = \lambda \psi,\quad \psi>0,\quad ||\psi||_{L^\infty}=1.
\end{equation}

The main results of the present work are stated as follows.

\begin{theorem}
	\label{thm:nonex}
	Let $f$ satisfy \ref{item:general}--\ref{item:regularity} for some $f_0$ and $a$, $c_0$ be the spreading speed of $f_0$, $\lambda$ be the supremum of the spectrum of $\partial_{xx}+a(x)$,
	and assume $\lambda > c_0^2$.
	
	\begin{enumerate}
		\item All entire solutions $u$ of \eqref{eq:main} with $0<u<1$, $\inf_{(t,x)\in \mR^2} u=0$ are bump-like. That is, there are $c,C_t>0$ such that $u(t,x)\le C_te^{-c|x|}$ for all $(t,x)\in \mR^2$.
		In particular, \eqref{eq:main} does not admit a transition front solution.
		
		\item Assume \ref{item:regularity} is replaced by the following: there exists $\zeta  \in (0,\taz)$ such that $f(x,u) = a(x)u$ for $u\in [0,\zeta]$.
		Then a nonzero bump-like solution of \eqref{eq:main} exists, is unique (up to a time-shift) among all solutions with $0<u<1$, $\inf _{(t,x)\in\mR^2}u=0$,
		and satisfies $u(t,x)=\zeta e^{\lambda t}\psi(x)$ for all $t\le 0$ after a time translation, where $\psi$ is the principal eigenfunction given in \eqref{eq:eigen}.
	\end{enumerate}
\end{theorem}

\begin{theorem}
	\label{thm:ex}
	Let $f,$ $c_0,$ and $\lambda$ be as in Theorem \ref{thm:nonex}.
	If $\lambda < c_0^2$, then \eqref{eq:main} admits a transition front, which is increasing-in-time and has a global mean speed $c_0$.
\end{theorem}

\begin{remark}
If \eqref{eq:main} does not admit a stationary solution other than $u\equiv 0, 1$ (i.e. $f(\cdot,\ta)\not\equiv 0$ for all $\ta \in (0,1)$), then any entire solution $u\not\equiv 1$ to \eqref{eq:main} must satisfy $\inf _{(t,x)\in \mR^2} u=0$ (see the proof of Lemma \ref{lem:limit}).
Hence in this case, Theorem \ref{thm:nonex} holds without the infimum assumption.
\end{remark}

\begin{remark}
It can been seen in the proof that Theorem \ref{thm:nonex} still holds if we replace $f(x,u)=f_0(u)$ in \ref{item:ignition} by $f(x,u)\le f_0(u)$ for $|x|\ge L$.
\end{remark}

\begin{remark}
The transition front $u$ constructed in Theorem \ref{thm:ex} satisfies $\sup _{t\in \mR} |X
(t) - c_0t|<\infty$, which is a stronger condition than \eqref{gms}. 
In fact, $u$ can be viewed as a perturbation of a traveling front $U$ of $f_0$.
By the asymptotic stability of ignition fronts, one can further show that $u$ converges to $U$ as $t\to\pm\infty$ in the sense that
\[ \lim_{t\to\pm \infty} \sup_{x\in\mR} | u(t,x) - U(x-c_0t -x_\pm) | = 0\quad \mbox{for some }x_\pm\in \mR. \]
\end{remark}

\begin{remark}	
In \cite{Zlatos}, infinitely many transition fronts with different global mean speeds were constructed in the case of Fisher-KPP reactions with weak inhomogeneity.
In our case, we only obtain a single transition front with a single global mean speed $c_0$ in Theorem \ref{thm:ex}. 
Indeed, by a comparison principle argument and stability of front speeds with respect to reactions, one can easily show that any transition front must have a global mean speed $c_0$. 
Though transition fronts in ignition media are expected to be unique \cite{Z12}, we do not know if the uniqueness still holds in our case of mixed ignition-monostable reactions. 

\end{remark}

Our method is structurally similar to that of \cite{Zlatos}.
The primary difference and technicality in this work is the lack of KPP structure in our reaction.
Indeed, KPP structure allows one to exploit the intimate connection between solutions of \eqref{eq:main} and those of the linearized equation $u_t=u_{xx}+a(x)u$ (which are super-solutions to \eqref{eq:main}).
Therefore, super- and sub-solutions found in \cite{Zlatos} are based on (generalized) eigenfunctions of the operator $\partial_{xx}+a(x)$ and exponential functions.
In our case, we clearly lack this convenience.
Hence, to overcome this, we base our super- and sub-solutions on the traveling front of the ignition reaction $f_0$ instead.
This modification introduces some difficulties in calculation, but we are still able to obtain a sharp result.

The body of the paper is organized as follows.
We prove the non-existence of transition fronts (Theorem \ref{thm:nonex}) in the coming section, and present the complimentary existence result (Theorem \ref{thm:ex}) in Section 3.

\section*{Acknowledgment}
All authors were supported in part by the NSF grant DMS-1056327, and thank Andrej Zlato\v{s} for directing this research.
CG gratefuly acknowledges the support of the Fannie and John Hertz Foundation.
CG, AM, and DW thank the Department of Mathematics at the University of Wisconsin--Madison for its hospitality during the 2013 REU ``Analysis and Differential Equations,'' where this research was mainly conducted.

\section{Non-existence for $\lambda >c_0^2$ (proof of Theorem \ref{thm:nonex})}
\label{sec:non}

As mentioned above, the methods of this section are based on those found in \cite{Zlatos}.
In particular, Theorem \ref{thm:nonex}, Lemmas \ref{lem:refined}, \ref{lem:new}, and their proofs are similar to Theorem 1.2, Lemmas 3.1, 3.2 \cite{Zlatos}.
The primary difference can be found in the proof of Lemma \ref{lem:new}.

Throughout this section, we assume $f,\ga,f_0, \taz,\tao, L, a$ are all from (F), and $\lambda > c_0^2$.
For $\eps\in (0,1)$, let $\lambda_\eps$ be the principal eigenvalue of the differential operator $\partial_{xx}+(1-\eps)a(x)$.
Since $\lim_{\eps\to 0^+} \lambda_\eps=\lambda>c_0^2$, we may fix $\eps>0$ such that $\lambda_\eps>c_0^2$, and
let $\zeta = \zeta(\eps)$ be given in \ref{item:regularity}.
For $M>0$, we let $\lambda_M = \lambda_{\eps,M}$ be the Dirichlet principal eigenvalue of $\partial_{xx}+(1-\eps)a(x)$ on $[-M,M]$, and $\psi_M\in C^2([-M,M])$ be the corresponding $L^\infty$-normalized eigenfunction:
\begin{gather}
	\lb{TS:psi}
	\psi_M'' + (1-\eps)a(x) \psi_M = \lambda_M \psi_M \mbox{ on }(-M,M),\\
	\psi_M(\pm M)=0,\quad ||\psi_M||_{L^\infty}=1. \nonumber
\end{gather}
Note that $\psi_M>0$ on $(-M,M)$ and $\lim_{M\to\infty} \lambda_M = \lambda_\eps>c_0^2$.
So we may again fix $M\ge L$ large so that $\lambda_M>c_0^2$. Finally, all constants involved depend on $c_0, M, \psi_M, \lambda_M, \zeta, \gamma, \taz$.

In the following, let $u\in (0,1)$ be an entire solution of \eqref{eq:main} with $\inf _{(t,x)\in \mR}u=0$.
In the proofs, we will frequently use the parabolic Harnack inequality for $u$.
Therefore, for $R,\sigma>0$, we let $k=k(R,\sigma)>0$ denote the Harnack constant such that
\begin{equation}
  \label{h-const}
	\min_{|x-x_0|\le R} u(t+\sigma,x) \ge k \max _{|x-x_0|\le R} u(t,x) \ge k u(t,x_0)
\end{equation}
holds for any $x_0\in \mR$.
We begin with the following simple fact.

\begin{lemma}
	\label{lem:limit}
	The solution $u$ satisfies $\lim_{t\to-\infty}u(t,x)=0$ locally uniformly.
\end{lemma}

\begin{proof}
	By the Harnack inequality, it suffices to show the limit for $x=0$.
	Assume the contrary, so that there exists $\al \in (0,1)$ and a sequence of times $\{t_n\}$ with $t_n\searrow -\infty$ such that $u(t_n,0)\geq \alpha$.
	Let $k=k(M,1)$ be the Harnack constant from \eqref{h-const},
	$\theta:=\min\{k\al, \zeta\}$, and extend the eigenfunction $\psi_M$ from \eqref{TS:psi} continuously to $\mR$ by setting $\psi_M\equiv 0 $ on $[-M,M]^c$.
	Since $||\psi_M||_{L^\infty}=1$, \eqref{h-const} (with $(R,\sigma,x_0,t)=(M,1,0,t_n)$) implies
	\begin{equation}
		\label{eq:Harnack_bd}
		u(t_n+1,x)\geq \theta\psi_M(x)\quad \mbox{for all }x\in \bbR.
	\end{equation}
	
	Now let $v:\mR^+\times \mR \to [0,1]$ be the solution to the Cauchy problem of \eqref{eq:main} with initial data $v(0,x) = \theta\psi_M(x)$.
	We claim that $v_t\ge 0$.
	By the comparison principle, it suffices to show $v(s,\cdot)\ge \theta\psi_M$ for all $s\ge 0$.
	Clearly this holds for all $x\in [-M,M]^c$ because $\psi_M\equiv 0$ in this region.
	If $x\in [-M,M]$ instead, observe that $w(t,x):=\theta\psi_M(x)$ is a (stationary) sub-solution of \eqref{eq:main} by (F3) and \eqref{TS:psi}. So the comparison principle shows that $v(s,x)\ge \theta\psi_M(x)$ for $(s,x)\in \mR^+\times [-M,M]$.
	This implies $v_t\ge 0$.
	Let $v_\infty(x):= \lim_{t\to\infty}v(t,x)$, which satisfies $v_\infty '' +f(x,v_\infty)=0$ on $\mR$ by parabolic regularity.
	Since $f\ge 0$, this forces $v_\infty \equiv \be$ for some constant $\be \in [\theta ,1]$.
	Now fix $s\in \mR$.
	By the comparison principle and \eqref{eq:Harnack_bd}, for all large $n$
	\begin{equation*}
		u(s,x)\ge v(s-t_n-1,x)\quad \mbox{for all }x\in \mR.
	\end{equation*}
	Letting $n\to\infty$, we find that $u(s,\cdot)\ge \be >0$ for all $s\in \mR$, contradicting $\inf_{(t,x)\in \mR^2} u=0$. Therefore, we must have $\lim_{t\to-\infty}u(t,0)=0$.
\end{proof}

With Lemma \ref{lem:limit}, after an appropriate time translation we may now assume
\begin{gather}
	\label{eq:zetabound}
	u(0,0)\leq \f \zeta 2 \psi_M(0).
\end{gather}
In the coming two lemmas, we establish some important bounds on $u$, which play a crucial role in the proof of Theorem \ref{thm:nonex}.

\begin{lemma}
	\label{lem:refined}
	For any $c\in (c_0, \sqrt{ \lambda_M} )$, there exists $C_0>0$ (independent of $u$) such that
	\begin{equation}
		\label{eq:refined_bound}
		u(t, x) \leq C_0 u(0,0) e^{c_0(x+ct)}, \quad \text{for }t\le -1,\; x\in [M,\sqrt{\laa_M}(-t-1)-M-1].
	\end{equation}
\end{lemma}

\begin{proof}
	Denote $u_0 := u(0,0)>0$, $\psi_0 := \psi_M(0)>0$, and $D\subset \mR^2$ the region described in \eqref{eq:refined_bound}.
	To show \eqref{eq:refined_bound}, we will prove the following estimate for some $C_0'>0$ (independent of $u)$:
	\begin{equation}
		\lb{TS:bound2.2}
		u(t, x) \leq C'_0 u_0 \sqrt{|t|}  e^{\sqrt{\lambda_M}(x+\sqrt{\lambda_M}t)},\quad \text{for }(t,x)\in D.
	\end{equation}
	One can easily show that \eqref{eq:refined_bound} follows from this with $C_0 := C'_0 \sup_{t \leq 0}\sqrt{|t|}e^{c_0(\sqrt{\lambda_M} - c)t}$ (which is finite because $\sqrt{\lambda_M}>c$).
	
	We prove \eqref{TS:bound2.2} by contradiction.
	Let $k=k(1,1)$ be the Harnack constant from \eqref{h-const}, and
	suppose there is $(t',x_0)\in D$ so that \eqref{TS:bound2.2} does not hold with  $C_0'$ given by
	\[ C_0':=\frac{\sqrt {4 \pi}}{k\psi_0}e^{\lambda_M+\sqrt{\lambda_M}(M+1)}. \]
	Let $t_0:=t'+1\le 0$ and 
	\begin{gather*}
		\be:=\f {x_0+M+1}{2|t_0|\sqrt{\lambda_M}},\quad 
		\eta:= C'_0ku_0\sqrt{|t'|}e^{\sqrt{\lambda_M}(x_0+\sqrt \lambda_Mt') }.
	\end{gather*}
	Observe that $\be \in (0,\f 12 ]$ as $(t',x_0)\in D$.
	Also, by \eqref{h-const} (with $(R,\sigma,t)=(1,1,t')$) and the opposite of \eqref{TS:bound2.2} we have $u(t_0,\cdot)\ge \eta \chi_{[x_0,x_0+1]}$.
	Applying the comparison principle ($u$ is a super-solution to the standard heat equation as $f\ge 0$), for all $x\in [-M,M]$ we have
	\begin{align}
		u(t_0+\be |t_0|,x) &\ge \f {\eta}{\sqrt{4\pi\be |t_0|}} \int _{x_0}^{x_0+1} e^{-\f{(x-y)^2}{4\be |t_0|}} dy\ge \f{\eta}{\sqrt{4\pi \be |t_0|}} e^{-\frac{(x_0+M+1)^2}{4\be|t_0|}} \nonumber \\
		&\ge {2u_0}{\psi_0}^{-1} e^{\sqrt{\lambda_M}(x_0+M+1+\sqrt{\lambda_M}t_0)-\frac{(x_0+M+1)^2}{4\be|t_0|}} \nonumber \\
		& = {2u_0}{\psi_0}^{-1} e^{\lambda_M(t_0+\be|t_0|)}.
		\lb{TS:2.2.3}
	\end{align}
	Here, the second inequality is due to $-M\le x \le y\le x_0+1$, for then $0\le y-x\le M+x_0+1$.
	Now let $v(t,x):=2u_0\psi_0^{-1}e^{\lambda _M t} \psi_M(x)$, which by \eqref{TS:psi} satisfies
	\begin{equation}\label{pdineq-23}
		v_t = v_{xx} + (1-\eps )a(x) v\quad \mbox{for all }(t,x)\in\mR\times (-M,M),\quad v(t,\pm M)=0.
	\end{equation}
	From \eqref{TS:2.2.3}, $||\psi_M||_{L^\infty}=1$, and \eqref{eq:zetabound}, we also have
	\begin{equation}\label{ineq-27}
		\begin{split}
			u(t_0+\be|t_0|,x)&\ge v(t_0+\be|t_0|,x), \\
			v(t,x)\le v(0,x)&=2u_0\psi_0^{-1}\psi_M(x)\le \zeta,
		\end{split}
	\end{equation}
	where the latter holds for all $t\le 0$ and $x\in \mR$.
	The latter with \eqref{pdineq-23} and \ref{item:regularity} shows that $v$ is a sub-solution of \eqref{eq:main} on $\mR^- \times (-M,M)$.
	Hence, \eqref{ineq-27} and the comparison principle (note that $t_0+\be |t_0|\le 0$ as $\be\in (0,\f 1 2]$ and $t_0\le 0$) yield
	\begin{equation*}
		u(0,x) \ge v(0 , x)=2u_0\psi_0^{-1}\psi_M(x) \quad  \mbox{for } x\in[-M,M].
	\end{equation*}
	Letting $x=0$ yields  the contradiction $u_0\ge 2u_0$ (as $u_0>0$).
	Therefore \eqref{TS:bound2.2} holds.
\end{proof}

With the estimate \eqref{eq:refined_bound}, we now further refine the bound for $u$ to show that it is bump-like for all large negative time.

\begin{lemma}
	\lb{lem:new}
	Under the same assumptions as Lemma \ref{lem:refined}, there exist $C>0$ and $\tau <0$ (both independent of $u$) such that
	\begin{equation}
		\lb{TS:new}
		u(t,x)\le C u(0,0) e^{-c_0|x|+c_0ct},\quad \mbox{for } t\le \tau,\, |x|\ge M.
	\end{equation}
\end{lemma}

\begin{remark}
We follow the argument in the proof of Theorem 1.2 in \cite{Zlatos}.
The fundamental difference lies in the definition of the super-solution $w=w_1+w_2^{s}$.
In \cite{Zlatos}, $w_2^{s}$ is chosen to be an exponential function (see the definition of $v_{t_0}$ in Section 3 \cite{Zlatos}).
In our case, $w_2^{s}$ is a (leftward-moving) traveling front of a perturbed reaction $f_\del$ (defined in the proof).
\end{remark}

\begin{proof}
	First of all, it suffices to prove \eqref{TS:new} for the case $x\ge M$.
	This is because $\til u(t,x):=u(t,-x)$ is still a global solution to \eqref{eq:main} with $\til f(x,u):=f(-x,u)$ in place of $f(x,u)$.
	Clearly, $\til f$ and $\til u$ satisfy (F) and \eqref{eq:zetabound} respectively.
	Applying \eqref{TS:new} to $\til u$ and $x\ge M$, we find \eqref{TS:new} for $u$ and $x\le -M$.
	
	Fix $c\in (c_0,\sqrt{\lambda_M})$.
	For $\del>0$ small, consider the perturbation of $f_0$ given by
	\begin{equation*}
		f_\delta(u) := \max_{v\in [u-\del,u+\del ]}f_0(v),
	\end{equation*}
	which is a Lipschitz ignition reaction with $\supp(f_\delta)= [\taz -\del,1+\del]$.
	Let $(U_\del,c_\del)$ be the traveling front solution to the following problem:
	\begin{gather}
		\lb{eq:tfeq2}
		U''_\del - c_\del U' + f_\del (U_\del ) = 0,\quad
		\lim_{y\to -\infty} U_{\delta}(y) = 0,\quad  \lim_{y\to \infty} U_{\delta}(y) = 1+\delta.
	\end{gather}
	Note that $U_\del$ is leftward moving, so $U_\del'>0$.
	By a simple argument using phase plane analysis and the uniqueness/stability of solutions to ODEs, one can easily show that $c_\del>c_0$ and $\lim_{\delta \searrow 0}c_{\delta}= c_0$.
	Thus we may fix $\del \in (0,\taz)$ so that $c_\del \in ( c_0, \sqrt {c_0c})$.
	Let $C_0$ be the constant from Lemma \ref{lem:refined} and $u_0:=u(0,0)$.
	Since $f_\del\equiv 0$ on $[0,\ta_0-\del]$ and \eqref{eq:tfeq2}, we can specify the translation of $U_\del$ so that
	\begin{equation}\label{Udel-trans}
		U_\del(x)=Ae^{c_\del x}\mbox{ whenever }U_\del(x)\le \ta_0-\del,\mbox{ where }A:=C_0u_0 e^{(c_0-c_\del )M}.
	\end{equation}
	For $s\in \mR$, define
	\begin{align*}
		w_1(t,x) := C_0u_0e^{-c_0(x-2M-ct)},\quad w_2^{s}(t,x) := U_\delta\rb {x + c_\delta t + \rb{\f {c c_0}{c_\delta} - c_\delta}s}.
	\end{align*}
	Let $\tau:=\min\{T_0,T_1\}$, where $T_0$, $T_1$ are given by
	\begin{equation*}
		T_0:= -\f {2M+1}{\sqrt {\laa_M}}-1 ,\qquad C_0 e^{c_0(M+cT_1)}=\del.
	\end{equation*}
	Here, $T_0$ is defined so that the interval from \eqref{eq:refined_bound} is non-empty for all $t\le T_0$.
	By Lemma \ref{lem:refined},
	\begin{equation}
		\lb{TS:comp1}
		u(t,M)\leq w_1(t,M)\quad \mbox{for }t\le \tau.
	\end{equation}
	We also claim that for all sufficiently negative $s\le 0$,
	\begin{equation}
		\lb{TS:comp2}
		u(s,x)\le w_2^{s}(s,x) \quad \mbox{for }x\in [M,\infty).
	\end{equation}
	We postpone the proof of this claim to first show \eqref{TS:new}.
	
	Let $w:=w_1+w_2^{s}$ and $D_{s} := [s,\tau]\times [M,\infty)$.
	Then $w$ is a super-solution to \eqref{eq:main} on $D_{s}$.
	After all, $w_1(t,x)\leq w_1(T_1,M)= u_0\delta\le \delta$ for all $(t,x) \in D_{s}$, so
	\begin{align*}
		w_t-w_{xx}&=\partial_t w_1-\partial_{xx}w_1+\partial_tw_2^{s}-\partial_{xx}w_2^{s}\nonumber \\
		&=c_0(c-c_0)w_1+f_\delta(w_2^{s})\geq f_\delta(w_2^{s})\nonumber \\
		&\ge f_0(w)=f(x,w).
	\end{align*}
	The last inequality follows from $0< w_1\leq \delta$ on $D_{s}$ and the definition of $f_\del$.
	By \eqref{TS:comp1}, \eqref{TS:comp2}, and  the comparison principle, we have $u\leq w_1+w^{s}_2$ on $D_{s}$, which holds for all large negative $s$.
	Observe that the argument of $U_\del$ in the definition of $w^{s}_2$ tends to $-\infty$ as $s\searrow -\infty$, since $c_\del <\sqrt {cc_0}$.
	Hence, $w^{s}_2\searrow 0$, and $u(t,x)\le w_1(t,x)$ for all $(t,x)\in (-\infty,\tau]\times [M,\infty)$.
	This is \eqref{TS:new} for $x\ge M$ if we set $C:=C_0e^{2c_0M}$.
	
	It remains to prove \eqref{TS:comp2}.
	Let $\xi_0 \in \mR$ satisfy $C_0u_0e^{c_0\xi_0}=\ta_0-\del$, and define
	\begin{equation*}
		W(s):= w_2^{s}(s,\xi_0-cs) = U_\del (\xi_0+c(c_0c_\del^{-1}-1)s),
	\end{equation*}
	which is continuous and satisfies $W(-\infty)=U_\del(\infty)=1+\del$ (as $c_0<c_\del$).
	Now, by this and $c<\sqrt{\lambda_M}$, we let $s_0\le -1$ so that for all $s\le s_0$ we have
	\begin{equation}\label{s0-def}
		W(s)\ge 1,\qquad \xi_0-cs\le -\sqrt {\lambda_M}-M-1-\sqrt{\lambda_M}s.
	\end{equation}
	Now we show that \eqref{TS:comp2} holds for all $s\le s_0$.
	Indeed, if $x\ge \xi_0-cs$, we have $w_2^{s}(s,x) \ge w_2^s(s,\xi_0-cs) =W(s)\ge 1 > u(s,x)$.
	Consider $x\in [M,\xi_0-cs)$. By \eqref{s0-def} and \eqref{eq:refined_bound}, it suffices to show that $w_2^s(s,x)\ge C_0u_0e^{c_0(x+cs)}$.
	Assume the contrary that it does not hold for some $x_0\in [M,\xi_0-cs)$.
	It then follows from the definition of $\xi_0$ that
	\begin{equation}\label{cont-05}
		w_2^s(s,x_0)< C_0u_0e^{c_0(x_0+cs)}\le C_0u_0e^{c_0\xi_0}=\ta_0-\del.
	\end{equation}
	On the other hand, by our translation for $U_\del$ from \eqref{Udel-trans}, $c_0<c_\del$ and $x_0\ge M$, it follows
	\begin{align*}
		w^{s}_2(s,x_0) &= U_\del(x_0+cc_0c^{-1}_\del s) =  Ae^{c_\del x_0+cc_0s}  \\
		&= C_0u_0 e^{c_0M+c_\del(x_0-M)+cc_0s}
		\ge C_0u_0e^{c_0(x_0+cs)}.
	\end{align*}
	This is a contradiction to the first inequality of \eqref{cont-05}. Hence, \eqref{TS:comp2} holds for all $s\le s_0$.
\end{proof}

We now prove our first main theorem.

\begin{proof}[Proof of Theorem \ref{thm:nonex}]
	(i) We have already shown in Lemma \ref{lem:new} that $u$ is bump-like for all $t\le \tau$.
	It follows from a standard argument using the comparison principle that $u$ remains bump-like for all time. This clearly implies that transition front does not exist.

	(ii) Denote $\varphi(t,x):=\zeta e^{\lambda t }\psi(x)$, which solves the following PDE globally in time:
	\begin{equation}
		\label{eq:linearized}
		\varphi_t = \varphi_{xx} + a(x)\varphi.
	\end{equation}
	To construct a bump-like solution, we let $u(t,x):= \varphi(t,x)$ for $t\le 0$, and propagate forward in time by solving the Cauchy problem for \eqref{eq:main} with initial data $u(0,x) = \zeta \psi(x)$ for $t >0$.
	Next we prove that this solution is unique up to time translation. Let $\til u\in (0,1)$ be another solution of \eqref{eq:main} with $\inf _{(t,x)\in \mR^2}\til u=0$.
	By Lemmas \ref{lem:limit} and  \ref{lem:new}, $\til u(t,\cdot)\to 0$ uniformly as $t\to-\infty$.
	Therefore, after a time-shift we may assume
	\begin{equation}\label{sup-trans}
		\sup _{(t,x)\in \mR^-\times \mR}\til u(t,x)\le \f \zeta 2 \psi_M(0).
	\end{equation}
	Let $\til \varphi(t,x):=\til u(t,x)$ for all $t\le 0$, and propagate forward in time as the solution of \eqref{eq:linearized} with $\til \varphi(0,x)=\til u(0,x)$.
	Since $f(x,u)=a(x)u$ for $u\in [0,\zeta]$ by the assumption, $\tilde \varphi$ is an entire solution of \eqref{eq:linearized}.
	By \cite[Proposition 2.5]{HP07}, we have $\til \varphi= q \varphi$ for some $q>0$, provided that Conditions (A), (H1), and (2.12) from \cite{HP07} are met (which will be shown shortly).
	Therefore, we have $\tilde \varphi(t,\cdot) = \varphi(t-T,\cdot)$ with $T:=-\lambda^{-1}\log({q}{\zeta}^{-1})$.
	Since $u\equiv \varphi$, $\tilde u \equiv \tilde \varphi$ for all $t\le 0$, it clearly follows $\tilde u(t,\cdot)=u(t-T,\cdot)$, which shows the uniqueness of solution.
	
	It remains to check all the conditions from \cite{HP07}.
	Note that (A) follows from $0\le a\le \ga\chi_{[-L,L]}$, and (H1) holds for the PDE \eqref{eq:linearized} because $\lambda >0$.
	To show \cite[(2.12)]{HP07}, we will prove that
	\begin{equation}
		\lb{TS:b1}
		\sup_{x\in\mR} \tilde \varphi(s,x)\le K \tilde \varphi(s,0) \quad \mbox{for all }s\in \mR,
	\end{equation}
	for some $K>0$ independent in time.
	Consider the above for $s\le 0$
	Let $\tilde \varphi^s(t,x):= \tilde \varphi(t+s,x)$, which again satisfies \eqref{eq:zetabound} (by \eqref{sup-trans}).
	Follow from Lemma \ref{lem:new}, $\til \varphi^s$ satisfies the estimate \eqref{TS:new}. With $t=\tau\le 0$, we find that
	\begin{equation*}
		\sup_{|x|\ge M}\til \varphi^s(\tau,x)\le C \til\varphi^s(0,0).
	\end{equation*}
	On the other hand, by the Harnack inequality \eqref{h-const}, we have $\max_{|x|\le M}\til \varphi^s(\tau,x)\le k^{-1}\til \varphi^s(0,0)$, where $k=k(M,-\tau)$.
	Hence, $\sup_{x\in \mR}\til \varphi^s(\tau,x)\le A\varphi^s(0,0)$ with $A:=\max\{C,k^{-1}\}$.
	Applying the comparison principle (noting that $w(t,x)=A\varphi^s(0,0)e^{\ga (t-\tau)}$ is a super-solution to \eqref{eq:linearized}), we find $\sup_{x\in \mR}\til\varphi^s (0,x)\le A e^{-\ga \tau} \til\varphi^s(0,0)$, which is \eqref{TS:b1} with $K:=Ae^{-\ga\tau}$.
	
	Now consider \eqref{TS:b1} with $s>0$.
	Decompose $\tilde \varphi(0,x)=\al\psi(x)+\psi^\perp(x)$, where $\psi$, $\psi^\perp$ are orthogonal in $L^2(\mR)$ (recalling that $\psi$ is the eigenfunction from \eqref{eq:eigen}).
	Let $\phi(t,x) := e^{-\lambda t}\til\varphi(t,x)$, which by \eqref{eq:linearized} satisfies
	\begin{equation*}
		\phi_t = (\partial_{xx}+a(x)-\lambda)\phi.
	\end{equation*}
	Since the principal eigenvalue $0$ of $\partial_{xx}+a(x)-\lambda$ is isolated, it is well-known that $\phi (t,\cdot)\to \al \psi $ uniformly as time progresses. This clearly implies \eqref{TS:b1} for $s>0$, as desired.
\end{proof}

\section{Existence for $\lambda<c_0^2$ (proof of Theorem \ref{thm:ex})}
\label{sec:ex}

In this section we will prove Theorem \ref{thm:ex} by finding a transition front solution of \eqref{eq:main} under the assumption $\lambda<c_0^2$.
The construction is based on the standard solution method mentioned in Section \ref{intro}---we first construct appropriate pairs of super- and sub-solutions, then obtain an entire solution $u$ sandwiched between them, and use their asymptotic behavior to show that $u$ is a transition front.

Throughout the section, $f,\ga,f_0, \taz,\tao, L$, and $a$ are all from Hypothesis (F),
and we extend $f(x,\cdot),f_0(\cdot)$ continuously to $\mR$ by setting $f_0(u)=f(x,u)\equiv 0$ for $u\notin [0,1]$.
For $\eps > 0$, let $\lambda^\eps>0,\psi^\eps> 0$ be the principal eigenvalue and (normalized) eigenfunction of the differential operator $\partial_{xx} + a(x) + 2\eps\chi_{[-L,L]}(x)$, which satisfy
\begin{gather*}
\partial_{xx}\psi^\eps+[a(x)+2\eps\chi_{[-L,L]}(x)] \psi ^\eps = \lambda^\eps \psi^\eps \mbox{ on }\mR\\
\lim_{|x|\to\infty}\epsi(x)=0,\quad ||\epsi||_{L^\infty}=1.
\end{gather*}
Since $a(x)=0$ for $|x|\ge L$ and $||\psi^\eps||_{L^\infty}=1$, $\psi^\eps$ satisfies the exponential bound
\begin{equation}
\lb{TS:psi-exp}
\epsi (x) \le \min\{1,e^{-\sqrt{\lambda^\eps}(|x|-L)}\}\quad \mbox{for all }x\in \mR.
\end{equation}
Note also that $\eps\mapsto \lambda^\eps$ is increasing and continuous, with $\lambda^0 = \lambda$ and $\lim_{\eps \to \infty}\lambda^\eps \to \infty$, so we may fix $\eps>0$ such that $\lambda^{\eps}\in (\f{c_0^2}{16},c_0^2)$, and let $\zeta = \zeta(\eps)\in (0,\taz)$ be from \ref{item:regularity}.
Let $U$ be the unique traveling front of $f_0$ in the sense of \eqref{tfeq} with $U(0)=\f {\taz} 2$.
Given $y\in \mR$, we define
\begin{equation}
\lb{TS:v1-def}
v^{y}(t,x) := U(x-c_0t + y),
\end{equation}
which satisfies \eqref{eq:main} with $f(x,u)$ replaced by $f_0(u)$.
Finally, let \begin{equation}
\lb{TS:om-eta-def}
\omega := \inf_{x\in[-L,L]} \psi^\eps(x),\quad \eta: = c_0\inf \left\{-U'(s);\;s\in \mR,\,U(s)\in \left[\frac{\taz}{2},1-\tao\right]\right\}.
\end{equation}
All the constants involved in this section will depend on $\ga, \taz,\tao, L, \eps, \zeta, U, c_0, \omega,$ and $\eta$.

We begin with the construction of sub- and super-solutions for $t\le 0$.
\begin{lemma}
	\label{lem:supersub}
	\begin{enumerate}
		\item For all $y\ge L$, $v^y$ given in \eqref{TS:v1-def} is a sub-solution to \eqref{eq:main} on $(-\infty,0)\times \mR$.
		\item  There exists $y_0\ge L+c_0\be(0)$, such that $w$ given as follows is a super-solution to \eqref{eq:main} on $(-\infty,0)\times \mR$\emph{:}
		\begin{gather*}
		w(t,x) := v^{y_0}(t+\be(t),x) + \phi(t,x),\\
		\beta(t):= \frac{\gamma\zeta}{4\sqrt{\ela}c_0\eta}e^{2\sqrt{\ela}c_0 t},\quad
		\phi(t,x):= \f \zeta 2 e^{\sqrt{\ela}c_0 t}\psi^\eps(x).
		\end{gather*}
	\end{enumerate}
\end{lemma}

\begin{proof}
	(i) Abbreviate $v=v^y$, recalling that $v$ satisfies \eqref{eq:main} with $f(x,u)\equiv f_0(u)$.
	Then we must show that
	\begin{equation}
	\label{eq:sub1}
	v_t-v_{xx}-f(x,v)=f_0(v)-f(x,v)\leq 0\quad \mbox{for }(t,x)\in (-\infty,0)\times \mR.
	\end{equation}
	For $x\in [-L,L]^c$, we have $f(x,v)\equiv f_0(v)$ by \ref{item:ignition}, and the above holds trivially.
	If $x\in [-L,L]$ instead, by $U'<0$ we have $v(t,x)=U(x-c_0t+y)\le U(y-L)\le U(0)= \f {\taz}2$.
	Therefore $f_0(v)=0$, and \eqref{eq:sub1} holds.
	
	(ii) Let $y_0\in \mR$ be the unique number such that
	\begin{equation}\label{y0-def}
	U(y_0- c_0\be(0)-L) = \min \cb{\f\zeta 2, \frac{\eps \omega \zeta }{2(\ga +\eps )}}.
	\end{equation}
	Clearly, $y_0\ge L+c_0\be(0)$ because $U(0)=\f {\taz}{2} >\frac{\zeta}{2}$ and $U'<0$.
	Abbreviate $v^{y_0} = v$. Then we must show for all $(t,x)\in (-\infty,0)\times \mR$ that
	\begin{align}
	w_t-w_{xx}-f(x,w)
	=&f_0(v) + v_t\,\beta' + \sqrt{\ela}\left(c_0-\sqrt{\ela}\right)\phi\nonumber \\
	& +[a(x)+2\eps\chi_{[-L,L]}(x)] \phi - f(x,w)\geq 0,
	\label{eq:super1}
	\end{align}
	where $v$ and $v_t$ are evaluated at $(t+\beta(t),x)$ in what follows.
	Consider first $x\in [-L,L]$.
	Since the first four terms in \eqref{eq:super1} are nonnegative, it suffices to show that $[a(x)+2\eps]\phi\geq f(x,w)$.
	Note that $\phi(t,x)\le \frac{\zeta}{2}$ as the eigenfunction $\psi^\eps$ satisfies $||\psi^\eps||_{L^\infty}=1$.
	Moreover, by $U'<0$, $\be '>0$ and \eqref{y0-def}, we have
	\begin{equation*}
	v(t+\be(t),x)\le U(y_0-c_0\be(0)-L)\le \frac{\zeta}{2}.
	\end{equation*}
	Therefore $w=v + \phi \le \zeta$, and $f(x,w)\leq\left[a(x)+\eps\right]w$ by \ref{item:regularity}.
	On the other hand, by \eqref{tfeq} and $f_0\equiv 0$ on $[0,\taz]$, we have $U(y)=\frac{\taz}{2}e^{-c_0y}$ for $y\ge 0$.
	\eqref{y0-def} then implies
	\begin{align*}
	v(t+\be(t),x)&\le U(y_0-c_0\be(0)-L-c_0t)\\
	&=U(y_0-c_0\be(0)-L)e^{c_0^2 t} \le \f {\eps \omega \zeta}{2(\ga+\eps )} e^{c_0^2t}.
	\end{align*}
	Combining this, $\sqrt{\lambda^\eps}<c_0$, $a(x)\le \ga$ and $\psi^\eps(x)\ge \omega $ on $[-L,L]$ (by \eqref{TS:om-eta-def}), we compute
	\begin{align*}
	[a(x)+2\eps]\phi-f(x,w)
	&\geq [a(x)+2\eps]\phi-\left[a(x)+\eps\right]w\\
	&= \eps\phi-\left[a(x)+\eps\right]v 
	\geq \f \zeta 2 \left[\eps \phi-\left[a(x)+\eps\right]\left(\frac{\eps\omega}{\ga+\eps }e^{c_0^2t}\right)\right]\\
	&\geq \f {\eps\zeta }2 e^{\sqrt{\lambda^\eps}c_0 t} [\psi^\eps(x)-\omega]\geq 0.
	\end{align*}
	Hence \eqref{eq:super1} holds for $(t,x)\in (-\infty,0)\times [-L,L]$.
	
	Now consider \eqref{eq:super1} for $x\in [-L,L]^c$.
	Since $f(x,w)=f_0(w)$ and the first four terms of \eqref{eq:super1} are non-negative, it suffices to show that
	\begin{align}
	f_0(v) + v_t\,\beta' - f_0(w)\geq 0.
	\label{eq:super2}
	\end{align}
	We consider three cases for the value of $v$.
	When $v\leq \frac{\taz}{2}$, we have $w=v+\phi\leq \taz$ (recalling that $\phi\le \frac{\zeta}{2}$), so $f_0(w)=0$ and \eqref{eq:super2} holds.
	If $v\geq 1-\tao$, $f_0(v)\geq f_0(w)$ because $f_0$ is non-increasing on $[1-\tao,\infty)$ by \ref{item:ignition}; \eqref{eq:super2} follows.
	Finally, suppose $\frac{\taz}{2} \leq v\leq 1-\tao$.
	In this case, by $y_0\ge c_0\be(0)+L$, $U'<0$, and $U(0)=\f {\taz}2$, we have $x\le c_0t-L$.
	Using $\mbox{Lip}(f)=\ga$, \eqref{TS:psi-exp}, and $v_t\ge \eta$ from \eqref{TS:om-eta-def}, we deduce
	\begin{align*}
	|f_0(v+\phi)-f_0(v)|&\leq \gamma \phi(t,x) \le \f {\ga\zeta } 2 e^{\sqrt{\ela}(x+L)+\sqrt{\ela}c_0 t}\le \f {\ga\zeta } {2} e^{2\sqrt{\ela}c_0t} \\
	&\leq v_t\frac{\gamma \zeta }{2\eta}e^{2\sqrt{\ela}c_0 t} =  v_t\,\beta'.
	\end{align*}
	Hence \eqref{eq:super2} always holds, and $w$ is a super-solution to \eqref{eq:main} on $(-\infty,0)\times\mR$.
\end{proof}

Next, we construct super- and sub-solutions for $t\ge 0$.

\begin{lemma}\lb{lem:sup-t+}
	\begin{enumerate}
		\item There exist $y_1\le -L$ and $B_1>0$ such that
		$\til w$ given as follows is a super-solution to \eqref{eq:main}
		on $(0,\infty)\times \mR$\emph{:}
		\begin{gather*}
		\til w(t,x):= v^{y_1} (t+\be_1(t),x)+\phi_1 (t,x),\\
		\be_1 (t):=B_1(1-e^{-c_0^2t/8}),\quad \phi_1(t,x) := e^{-\f{c_0}{4} (x-L -\f{c_0}{2} t)}.
		\end{gather*}
		\item For all $y\in \mR$, there exists $B_2=B_2(y)>0$ such that
		$\til v^y$ given as follows is a sub-solution to \eqref{eq:main}
		on $(0,\infty)\times \mR$\emph{:}
		\begin{gather*}
		\tilde v^y(t,x):=  v^{y} (t+\be_2(t),x) -\phi_2 (t,x),\\
		\be_2 (t):=B_2 e^{-c_0^2t/8},\quad \phi_2(t,x):= \f{16\ga }{c_0^2} e^{-\f{c_0}{4} (x-L-\f{c_0}{2}t)}.
		\end{gather*}
	\end{enumerate}

\end{lemma}

\begin{proof}
	(i) Let $y_1,\ell,B_1\in \mR$ be defined as follows:
	\begin{equation}
	\lb{TS:3.2-def}
	\phi_1(0,-y_1)=\f{\taz}{2},\quad U(\ell)=1-\tao,\quad B_1 := \f{8\ga }{c_0^2 \eta }e^{-\f{c_0}{4}(\ell-y_1 -L)}.
	\end{equation}
	Note that $y_1\le -L$ because $\phi_1(0,L)=1$. As before, we abbreviate $v=v^{y_1}$, and we need to show for all $(t,x)\in(-\infty,0)\times \mR$ that
	\begin{equation}
	\label{eq:tgreater}
	\til w_t - \til w_{xx} - f(x,\til w) =  v_t\,\be_1' + \frac{c_0^2}{16}\phi_1+f_0(v)- f(x,\til w)\geq 0,
	\end{equation}
	where $v$ and $v_t$ are evaluated at $(t+\beta_1(t),x)$ for the remainder of this part.
	Note that the first three terms are nonnegative.
	If $x\le L$, then \eqref{eq:tgreater} holds by $f(x,\til w) =0$ as $\til w(t,x)\geq \phi_1(t,x)\ge 1$.
	Now consider $x>L$, noting that $f(x,\til w) = f_0(\til w)$.
	We again consider three cases for the value of $v$.
	When $v\geq 1-\tao$, we have $f_0(v)-f_0(\til w)\geq 0$ because $f_0$ is non-increasing on $\left[1-\tao,\infty\right)$.
	\eqref{eq:tgreater} follows.
	If $v\leq \frac{\taz}{2}$, then \eqref{TS:v1-def} and $U(0)=\f {\taz}{2}$ imply $x-c_0t \ge -y_1$.
	Therefore,
	\begin{equation*}
	\phi_1(t,x)\le \phi_1(2t,x)=\phi_1(0,x-c_0t)\le \phi_1(0,-y_1)=\frac{\taz}{2}.
	\end{equation*}
	It follows that $\til w= v+\phi_1 \le \taz$, so $f_0(\til w)=0$, and \eqref{eq:tgreater} holds again.
	Finally, suppose $v\in [\frac{\taz}{2}, 1-\tao] $.
	From \eqref{TS:3.2-def}, we again have $x-c_0t \geq \ell-y_1 $.
	Using the definition of $\eta$ from \eqref{TS:om-eta-def},
	\begin{align*}
	|f_0(v)-f_0(\til w)|&\leq \ga \phi_1(t,x) = \ga e^{-\f{c_0}{4} (x-L -\f{c_0}{2} t)}\\
	&\leq \gamma e^{-\f{c_0 ^2}{8}t-\f{c_0}{4}(\ell-y_1 -L)} \le v_t \f{\ga}{\eta} e^{-\f{c^2_0}{8}t-\f{c_0}{4}(\ell-y_1 -L)} = v_t\,\beta_1'.
	\end{align*}
	\eqref{eq:tgreater} again follows, so $\til w$ is a super-solution to \eqref{eq:main} on $(0,\infty)\times \mR$.
	
	(ii) Define $\til \ta,\til \eta,\til \ell, B_2=B_2(y)$ as follows:
	\begin{align*}
	\til \theta := \min\cb{\f{\taz}{2},\f{\tao}{2},\f{c_0^2 \tao}{32\ga}},&\quad \til \eta:= c_0\inf\left\{-U'(s):s\in \mR,\,U(s)\in \left[\tta,1-\tta\right]\right\},\\
	U(\tl)= 1-\tta,&\quad \tB:= \f{2^7 \ga^2}{c_0^4\til \eta }e^{-\f{c_0}{4}\left(\tl-y -L\right)}.
	\end{align*}
	Again abbreviate $v=v^{y}$, $\til v=\til v^y$.
	We will show for all $(t,x)\in (0,\infty)\times \mR$ that
	\begin{align}
	\til v_t- \til v_{xx} - f(x, \til v ) = v_t\,\tbe_2' - \frac{c_0^2}{16}\tphi_2  + f_0(v) - f(x, \til v)\leq 0,
	\label{eq:subsol}
	\end{align}
	where $v$ and $v_t$ are evaluated at $(t+\tbe_2(t),x)$ for the rest of the proof.
	Note that $f_0(v)$ is the only term in \eqref{eq:subsol} which is potentially positive.
	Consider $x\in [-L,L]$, where we have $\tphi_2(t,x)\geq \tphi_2(0,L) = {16\ga}{c_0}^{-2}$.
	Then \eqref{eq:subsol} follows from
	\begin{equation*}
	- \frac{c_0^2}{16}\tphi_2 + f_0(v) \le -\ga+f_0(v) \leq 0.
	\end{equation*}
	Now consider $x\in[-L,L]^c$, where $f(x,\til v) = f_0(\til v)$.
	\eqref{eq:subsol} holds whenever $v \leq \tta(\le \taz)$, as $f_0(v)=0$.
	If $v\in[\tta, 1 - \tta]$, then $x\ge c_0 t+\tl-y$ because $U(\til \ell)=1-\til \ta$. It then follows from $v_t \ge \til \eta>0 $ that
	\begin{align*}
	\abs{f_0(v) - f_0(\til v)} &\leq \gamma \tphi_2(t,x) =  \f{16\ga^2}{c_0^2}e^{-\f{c_0}{4}(x-L-\f{c_0}{2}t)}\le \f{16\ga^2}{c_0^2} e^{-\f{c_0^2}{8}t-\f{c_0}{4}(\tl-y-L)} \\
	&\le v_t\f{16\ga^2}{c_0^2\til \eta}e^{-\f{c_0^2}{8}t-\f{c_0}{4}(\tl-y-L)} = -v_t\,\tbe_2',
	\end{align*}
	which implies \eqref{eq:subsol}.
	Finally, consider $v \geq 1 - \tta(\ge 1-\f{\tao}{2})$.
	If $\tphi_2 \leq \frac{\tao}{2}$ then $\til v \geq  1 - \tao$.
	Since $f_0$ is non-increasing on $[1 - \tao,\infty)$, $f_0( v) \leq f_0( \til v)$, and \eqref{eq:subsol} holds again.
	If $\tphi_2\geq \f{\tao}{2}$, then
	\begin{equation*}
	\frac{c_0^2}{16}\tphi_2\geq \frac{c_0^2\tao}{32}\geq \gamma \tta \geq \gamma \abs{v - 1} \geq f_0(v).
	\end{equation*}
	Hence \eqref{eq:subsol} holds in all cases, as claimed.
\end{proof}

With the requisite super- and sub-solutions in place, we may now establish Theorem \ref{thm:ex}.

\begin{proof}[Proof of Theorem \ref{thm:ex}]
	In what follows, $y_0,y_1, B_1, \be, \phi, w,\be_1, \phi_1,\til w,$ and $\phi_2$ are constants and functions defined in Lemmas \ref{lem:supersub}(ii) and \ref{lem:sup-t+}, and let $v=v^{y_0}$ be given as in \eqref{TS:v1-def}.
	
	The construction of an entire solution $u$ follows the standard procedure.
	For $n\in \mathbb N$, let $u_n$ be the solution to \eqref{eq:main} on $(-n,\infty)\times \mR$ with initial data $u_n(-n,x) = v(-n,x)$.
	First, observe that $u_n$ is increasing in time.
	Indeed, since $v$ is a sub-solution of \eqref{eq:main} on $(-\infty,0)\times \mR$ with $v_t>0$ by Lemma \ref{lem:supersub}, the comparison principle implies $u_n(t,\cdot)> u_n(0,\cdot)$ for all $t> 0$.
	It follows thast $\partial_t u_n> 0$ by the maximum principle.
	Moreover, since $v(-n,\cdot)\le w(-n,\cdot)$ by their definition, Lemma \ref{lem:supersub} and the comparison principle ensure that
	\begin{equation*}
	v(t,\cdot)\le u_n(t,\cdot)\le w(t,\cdot)\quad\mbox{for all }t\in [-n,0].
	\end{equation*}
	By parabolic regularity, we obtain an increasing in time entire solution $u$ to \eqref{eq:main} as a locally uniform limit along a subsequence of $(u_n)$ satisfying
	\begin{equation}
	\lb{TS:ex-comp}
	v(t,\cdot)\le u(t,\cdot)\le w(t,\cdot)\quad \mbox{for all }t\le 0.
	\end{equation}
	
	Next, we check that $u$ fulfills \eqref{item:trans_lims}.
	It obviously holds for $t\le 0$ by \eqref{TS:ex-comp} and the limit behavior of $v,w$ at $\pm \infty$.
	For $t>0$, the first limit of \eqref{item:trans_lims} still holds because $u_t>0$ and $u<1$.
	To prove the second limit condition, we first claim that
	\begin{equation}
	\lb{TS:ex-comp1}
	u(0,x)\le \min\{ 1, w(0,x) \} \le \til w(0,x).
	\end{equation}
	From this, Lemma \ref{lem:sup-t+}(i), and the comparison principle,
	\begin{equation}
	\lb{TS:ex-comp3}
	u(t,\cdot)\le \til w(t,\cdot)\quad \mbox{for all }t\ge 0.
	\end{equation}
	The second limit of \eqref{item:trans_lims} then follows from $\lim_{x\to\infty}\til w(t,x) =0$ and $u>0$.
	Now consider \eqref{TS:ex-comp1}.
	The first inequality is simply \eqref{TS:ex-comp} with $t=0$.
	For the second, when $x\le L$, $\til w(0,x)\ge \phi_1(0,x) \ge 1$.
	For $x>L$, \eqref{TS:psi-exp}, $\lambda^\eps >\f{c_0^2}{16}$, $U'<0$, and $y_1\le y_0-c_0\be(0)$ (by Lemmas \ref{lem:supersub}(ii) and \ref{lem:sup-t+}(i)) imply
	\begin{align*}
	w(0,x) &= U(x+y_0-c_0\beta(0)) + \frac{\zeta}{2}\psi^\eps(x)\le U(x+y_1) + e^{-\f{c_0}{4}(x-L)} = \til w(0,x).
	\end{align*}
	Therefore \eqref{TS:ex-comp1} holds.
	This completes the proof of \eqref{item:trans_lims}.
	
	It remains to show the bounded front width condition \eqref{item:trans_width}.
	Fix $\mu \in (0,\f 12 )$ and let
	\begin{equation*}
	X^-(t):=\inf \{x\in \mR:u(t,x)\le 1-\mu\},\quad X^+(t):=\sup \{x\in \mR:u(t,x)\ge \mu\},
	\end{equation*}
	for then we have $L_\mu(t) = X^+(t)-X^-(t)$.
	We will show that $L_\mu(t)$ is uniformly bounded in $t\in \mR$ by considering these three cases: $t\le 0$, $t> t_\mu,$ and $t\in (0,t_\mu]$, where $t_\mu$ will be defined shortly.
	For the first case $t<0$, let
	\begin{gather*}
	\rho_-:= U^{-1}(1-\mu),\quad 
	 \rho_+:= \max\left\{ U^{-1}\left (\f \mu 2\right), \f 1{\sqrt{\ela}}\left |\log \f \mu 2\right|+L+y_0-c_0\be(0) \right\}.
	\end{gather*}
	For all $x< c_0 t + \rho_--y_0$, by \eqref{TS:ex-comp} we have
	\begin{equation*}
	u(t,x)\ge v(t,x) = U(x-c_0t +y_0)> U(\rho_-)= 1-\mu.
	\end{equation*}
	Therefore $X^-(t) \ge c_0t +\rho_--y_0$.
	Now, suppose $x> c_0t+\rho_+ -y_0 +c_0\be(0)$.
	Combining \eqref{TS:psi-exp}, $||\psi^\eps||_{L^\infty}=1$, \eqref{TS:ex-comp} and the definition of $\rho_+$, we compute
	\begin{align*}
	u(t,x) &\le w(t,x) = U(x-c_0(t+\be(t))+y_0)+\f \zeta 2 e^{c_0\sqrt{\ela}t}\epsi(x) \\
	&\le  U(x-c_0(t+\be(0))+y_0)+ \f \zeta 2e^{-\sqrt{\ela}(x-L-c_0t)} \\
	&< U(\rho_+)+ e^{-\sqrt{\ela}(\rho_+-L-y_0+c_0\be(0))}\le \mu .
	\end{align*}
	This implies $X^+(t) \le c_0t +\rho_+ - y_0+c_0\be(0)$, so
	\begin{equation}
	\lb{TS:L-1}
	L_\mu(t) \le c_0\be(0)+\rho_+ - \rho_-\quad \mbox{for }t\le 0.
	\end{equation}
	
	We now define $t_\mu$.
	Given $y\in \mR$, let $\bar v^y(t,x):= v^{y}(t,x)- \phi_2(t,x)$, where $v^y,\phi_2$ are from \eqref{TS:v1-def} and Lemma \ref{lem:sup-t+}(ii).
	One can easily verify that $M(y):=\sup_{x\in \mR} \bar v^{y}(0,x)$ is continuous, non-increasing in $y\in \mR$, $\lim_{y\to-\infty} M(y)=1$, $\lim_{y\to\infty}M(y)=0$, and the supremum is achieved somewhere if $M(y)>0$.
	So we may fix $y_2=y_2(\mu)$ such that $M(y_2)=1-\mu$.
	For the remainder of the proof, we abbreviate $\bar v = \bar v ^{y_2}$ and let $B_2=B_2(y_2)$, $\tbe_2$, and $\til v=\til v^{y_2}$ be from Lemma \ref{lem:sup-t+}(ii).
	Let $x_\mu\in \mR$ be a maximizer so that $\bar v (0,x_\mu)=1-\mu$, and define
	\begin{equation}\label{t-mu-def}
	t_\mu := \inf \cb{t>0: u(t,\cdot)\ge \max \{\til v(0,\cdot), (1-\mu)\chi_{(-\infty,x_\mu]}\}  }\ge 0.
	\end{equation}
	We claim that $t_\mu$ is finite.
	After all, we have $\til v(0,\cdot)\le (1-\mu')\chi_{I}$ for some bounded interval $I\subset \mR$ and $\mu'\in (0,1)$.
	Recall that $u_t>0$ and the limit condition \eqref{item:trans_lims} holds.
	Therefore, we have $u(t,\cdot)\nearrow 1$ uniformly on each $(-\infty,R)$, $R\in \mR$, which implies $t_\mu<\infty$.
	
	Now consider the front width for $t>t_\mu$.
	Combining \eqref{t-mu-def}, Lemma \ref{lem:sup-t+}(ii), the comparison principle, and $\tbe_2>0$, we have
	\begin{equation}
	\lb{TS:ex-comp2}
	u(t_\mu+t,\cdot)\ge \til v(t,\cdot)\ge \bar v(t,\cdot)\quad \mbox{for all }t\ge 0.
	\end{equation}
	Recall $B_1$ from Lemma \ref{lem:sup-t+}(i) and set
	\begin{equation*}
	\til\rho_+ := \max \cb{U^{-1}\rb{\f \mu 2}, \f{4}{c_0} \abs{\log\f{\mu}{2}} +y_1+L-c_0B_1}.
	\end{equation*}
	Then $X^+(t) \le c_0(t+B_1) +\til \rho_+-y_1$.
	Indeed, for $x> c_0(t+B_1)+\til\rho_+-y_1$, \eqref{TS:ex-comp3} and $\be_1< B_1$ imply
	\begin{align*}
	u(t,x)&\le \til w(t,x) = U(x-c_0(t+\be_1(t))+y_1) + e^{-\f{c_0}{4} (x-L-\f{c_0}{2}t)}\\
	&< U(\til\rho_+)+ e^{-\f{c_0}{4}(c_0B_1+\til \rho_+-y_1-L)} \le \mu.
	\end{align*}
	On the other hand, we claim that $X^-(t)\ge c_0(t-t_\mu)+x_\mu$, implying:
	\begin{equation}
	\lb{TS:L-2}
	L_\mu(t) \le c_0(t_\mu+ B_1) +\til \rho_+ -y_1 - x_\mu,\quad \text{for }t\in (t_\mu,\infty).
	\end{equation}
	To prove the claimed bound, it suffices to check that $u(t,x)\ge 1-\mu$ for all $x< c_0(t-t_\mu)+x_\mu$.
	If $x\le x_\mu$, then $u_t>0$ and \eqref{t-mu-def} show that $u(t,x)> u(t_\mu ,x)\ge 1-\mu$.
	Now consider $x\in (x_\mu, c_0(t-t_\mu)+x_\mu)$.
	Note that $v^{y_2}(t,x_\mu+c_0t)=U(x_\mu+y_2)$, while $t\mapsto \tphi_2(t,x_\mu+c_0t)$ is decreasing.
	Hence, their difference $\bar v (t,x_\mu+c_0t) $ is increasing in $t$, and
	\begin{equation*}
	\bar v(t,x_\mu+c_0t)> \bar v(0,x_\mu)=1-\mu\quad  \mbox{ for all }t>0.
	\end{equation*}
	Let $t_*:=c_0^{-1}(x-x_\mu)\in (0,t-t_\mu)$.
	Then by $u_t> 0$ and \eqref{TS:ex-comp2}, it follows that
	\begin{equation*}
	u(t,x) > u(t_\mu+t_*,x)=u(t_\mu+t_*,x_\mu+c_0t_*)\ge \bar v(t_*,x_\mu+c_0t_*) > 1-\mu.
	\end{equation*}
	This proves the claim.
	
	Finally, consider $t\in (0,t_\mu]$.
	Since $u_t>0$, the width is bounded by
	\begin{equation*}
	L_\mu(t)\le X^+(t_\mu)-X^-(0) \le c_0 (t_\mu+B_1) +\til\rho_+-\rho_- +y_0-y_1.
	\end{equation*}
	With this, \eqref{TS:L-1}, and \eqref{TS:L-2}, $L_\mu(t)$ is uniformly bounded for all $t\in \mR$.
	This concludes the proof of \eqref{item:trans_width}.
	Therefore $u$ is an increasing-in-time transition front solution of \eqref{eq:main}. 
	It also obviously holds from the comparisons \eqref{TS:ex-comp}, \eqref{TS:ex-comp3} and \eqref{TS:ex-comp2} that $u$ has a global mean speed $c_0$. 
	This completes the proof of Theorem \ref{thm:ex}. 
\end{proof}

\bibliographystyle{amsplain}
\bibliography{ign_mon_front_arxiv.bib}

\end{document}